\documentclass[12pt]{amsart}
\usepackage{latexsym}
\usepackage{amsmath}
\usepackage{amsfonts}
\usepackage{amssymb}
\usepackage{amsthm}
\theoremstyle{plain}
\newtheorem{theorem}{Theorem}[section]
\newtheorem{lemma}[theorem]{Lemma}

\newtheorem*{theorem*}{Theorem}
\newtheorem*{Correspondence1}{Furstenberg's Correspondence Principle }

\newtheorem*{PSzemeredi}{Polynomial Szemer\'edi Theorem}

\theoremstyle{definition}
\newtheorem{definition}[theorem]{Definition}
\newtheorem*{example}{Example}

\theoremstyle{remark}

\newcommand{\E}{\mathbb{E}}
\newcommand{\bbE}{\mathbb{E}}
\newcommand{\N}{\mathbb{N}}
\newcommand{\bbN}{\mathbb{N}}

\newcommand{\bbT}{\mathbb{T}}
\newcommand{\bbZ}{\mathbb{Z}}

\newcommand{\cB}{\mathcal{X}}

\newcommand{\cD}{\mathcal{D}}

\newcommand{\cK}{\mathcal{K}}

\newcommand{\cX}{\mathcal{X}}

\newcommand{\cY}{\mathcal{Y}}

\begin{document}
\title{Ergodic averages for independent polynomials and applications}
\author{Nikos Frantzikinakis and Bryna Kra}
\begin{abstract}

Szemer\'edi's Theorem states that a set of integers with positive
upper density contains arbitrarily long arithmetic progressions.
Bergelson and Leibman generalized this, showing that sets of
integers with positive upper density contain arbitrarily long
polynomial configurations;  Szemer\'edi's Theorem corresponds to
the linear case of the polynomial theorem. We focus on the case
farthest from the linear case, that of rationally independent
polynomials.  We derive results in ergodic theory and in
combinatorics for rationally independent polynomials, showing that
their behavior differs sharply from the general situation.
\end{abstract}

\address{Department of Mathematics,
Institute for Advanced Study, 1 Einstein Drive, Princeton, NJ
08540}
\address{Department of Mathematics, Northwestern University,
2033 Sheridan Road,  Evanston, IL 60208-2730}
\email{nikos@math.ias.edu} \email{kra@math.northwestern.edu}
\thanks{The  first author was partially supported by NSF grant
DMS-0111298 and the second author  by NSF grant DMS-0244994.}



\maketitle
\section{Introduction and results in ergodic theory}

\subsection{Background}
The celebrated theorem of Szemer\'edi~\cite{Sz} states that a
subset of the integers with  positive upper density\footnote{If
$\Lambda\subset\mathbb{N}$ we define the \emph{upper density}
$\bar{d}(\Lambda)=\limsup_{N\to\infty}|\Lambda\cap[1,N]|/N$.}
contains arbitrarily long arithmetic progressions.
Furstenberg~\cite{Fu1} drew the deep connection between
combinatorial questions and ergodic theory, showing that
Szemer\'edi's Theorem follows from an ergodic theorem, now known
as the multiple recurrence theorem.

A natural question is to find other configurations that must occur
in subsets of the integers with positive upper density.
Furstenberg~\cite{Fu2} and S\'ark\"ozy~\cite{Sa} independently
proved that if $\Lambda\subset \N$ has positive upper density and
$p(n)$ is an \emph{integer polynomial}, meaning it takes integer
values on the integers,  and if $p(0)=0$, then there exist
$x,y\in\Lambda$ such that $x-y=p(n)$ for some $n\in\N$. Bergelson
and Leibman established a far reaching generalization of this
result. They showed that if $\Lambda\subset \N$ has positive upper
density and $p_1,\ldots,p_k$ are integer polynomials with
$p_i(0)=0$ for $i=1,\ldots,k$, then there exists $n\in\bbN$ such
that
\begin{equation}\label{E:density}
\bar{d}\bigl(\Lambda\cap (\Lambda+p_1(n))\cap\cdots\cap
(\Lambda+p_k(n))\bigr)>0 \ .
\end{equation}

As with Furstenberg's proof of Szemer\'edi's Theorem, the
Polynomial Szemer\'edi Theorem follows from an ergodic theorem:
\begin{PSzemeredi}[Bergelson and Leibman~\cite{BL}]\label{T:PSzemeredi}
Let $(X,\cB,\mu,T)$ be  an invertible  measure preserving system and let $p_1$,
\dots, $p_k$ be integer polynomials with $p_i(0) = 0$ for $i = 1,
\ldots, k$. If $A\in\cB$ with $\mu(A)>0$, then
\begin{equation}
\label{eq:blerg}
\liminf_{N\to\infty}\frac{1}{N}\sum_{n=0}^{N-1}\mu\bigl(A \cap
T^{p_1(n)}A \cap\ldots\cap T^{p_k(n)}A \bigr) > 0 \ .
\end{equation}
\end{PSzemeredi}

Szemer\'edi's Theorem (and the ergodic theoretic proof by
Furstenberg) corresponds to the case that all the polynomials are
linear. We focus on the opposite case of \emph{rationally
independent} integer polynomials, meaning a set of integer
polynomials such that every nontrivial integer combination of the
polynomials is not constant. In some sense, this case is typical,
since a generic family of integer polynomials is rationally
independent. A particular example is any set of polynomials with
pairwise distinct degrees. We prove several results, some ergodic
and some combinatorial, for families of rationally independent
integer polynomials, focusing on the difference between this case
and that of a family of linear integer polynomials.

\subsection{Ergodic Results}
Studying the limiting behavior of the multiple ergodic averages
associated with \eqref{eq:blerg}  has been a central topic in
ergodic theory. Very recently, using methods from \cite{HK}
convergence was established for totally ergodic systems
in~\cite{HK4} and for general systems  in~\cite{Lei3} .
 The basic approach is to
 find an appropriate factor system, called a
\emph{characteristic factor}, that controls the limiting behavior
as $N-M \to\infty$ in $L^2(\mu)$ of the averages
\begin{equation}\label{E:averages}
\frac{1}{N-M}\sum_{n=M}^{N-1} T^{p_1(n)}f_1 \cdot \ldots\cdot
T^{p_k(n)}f_k \ .
\end{equation}
A characteristic factor is a factor such that the limit of the
averages remains unchanged when each function is replaced by its
projection on this factor.  The next step is to obtain a concrete
description for some well chosen characteristic factor in order to
prove convergence. For general polynomials, such a characteristic
factor can be described as an inverse limit of nilsystems (defined
in Section~\ref{S:preliminaries}). We show that characteristic
factors for rationally independent  integer  polynomials have a
significantly simpler structure.  In particular, in
Section~\ref{sec:char} we show that a characteristic factor for
rationally independent polynomials can be chosen to be an inverse
limit of rotations on finite abelian groups:
\begin{theorem}\label{T:characteristic} Let $(X, \mathcal{X}, \mu, T)$ be
an ergodic invertible  measure preserving system and
$p_1,\ldots,p_k$ be rationally independent integer polynomials.
Then the rational Kronecker factor $\mathcal{K}_{rat}$ (defined in
Section~\ref{S:preliminaries}) is a characteristic factor for the
$L^2(\mu)$-convergence of the averages \eqref{E:averages}, meaning
that if $f_1, \ldots, f_k \in L^\infty(\mu)$, the difference
\begin{equation}\label{E:1}
\frac{1}{N-M}\sum_{n=M}^{N-1} T^{p_1(n)}f_1\cdot\ldots\cdot
T^{p_k(n)}f_k- \frac{1}{N-M}\sum_{n=M}^{N-1}
T^{p_1(n)}\tilde{f}_1\cdot\ldots\cdot T^{p_k(n)}\tilde{f}_k\ ,
\end{equation}
where $\tilde{f}_i=\bbE(f_i\,|\,\mathcal{K}_{rat})$,
$i=1,\ldots,k$, converges to $0$ in $L^2(\mu)$ as $N-M\to\infty$.
\end{theorem}

For a given measure preserving  system $(X,\mathcal{B},\mu,T)$ and
functions $f_0, f_1,\ldots, f_k\in L^\infty(\mu)$, it was shown in
\cite{BHK} that the multicorrelation sequence
$$
a_n=\int f_0 \cdot T^nf_1\cdot\ldots\cdot T^{kn}f_k\ d\mu
$$
can be decomposed as a sum of a  $k$-step nilsequence and a
sequence that converges to zero in uniform density (all notions
defined in Section~\ref{SS:nil}).  We note that the original
statement in~\cite{BHK} is for $f_0 = f_1 = \ldots = f_k$, but the
same proof holds for different functions. Using
Theorem~\ref{T:characteristic} we prove an analogous result for
the multicorrelation sequence of independent polynomial iterates.
Moreover, in Section~\ref{SS:nil} we show that a significantly simpler class of
nilsequences suffices for the decomposition:
\begin{theorem}\label{T:nilsequence}
Let $(X,\cX,\mu,T)$ be an invertible  ergodic measure preserving
system and let $p_1,\ldots,p_k$ be rationally independent integer
polynomials with highest degree $d$. If $f_0$, $f_1$, $\ldots,$
$f_k$ $\in L^\infty(\mu)$, $n\in\mathbb{N}$ and
$$
a_n=\int f_0\cdot T^{p_1(n)}f_1\cdot \ldots \cdot T^{p_k(n)}f_k \
d\mu \ ,
$$
then $\{a_n\}_{n\in\bbN}$ is the sum of a $d$-step nilsequence and
a sequence that converges to zero in uniform density. Moreover,
the $d$-step nilsequence can be chosen to be of the form
$b_n=\phi(S^ne)$, where $S\colon G^d\to G^d$ is a unipotent affine
transformation, $G$ is a compact abelian group, $\phi\colon
G^d\to\mathbb{C}$ is continuous, and $e$ is the identity element
of $G^d$.
\end{theorem}

We also use Theorem~\ref{T:characteristic} to prove a multiple
recurrence result. We show that for a family of rationally
independent integer polynomials, the measure of the intersection
in \eqref{eq:blerg} is as large as possible ``frequently.'' More
precisely, a set $\Lambda\subset\N$ is \emph{syndetic} if there
exists $M\in\bbN$ such that every interval of length greater than
$M$  intersects $\Lambda$ nontrivially.  In
Section~\ref{sec:multrec} we show:

\begin{theorem}\label{C:T-2}
Let $(X,\cX,\mu,T)$ be an invertible  measure preserving system,
$p_1,\ldots,p_k$ be rationally independent integer polynomials
with $p_i(0)=0$ for  $i=1,\ldots,k$, and  $A\in\cB$. Then  for
every $\varepsilon>0$, the set
$$
\left\{n\in\N\colon \quad \mu\left( A\cap
T^{p_1(n)}A\cap\cdots\cap T^{p_k(n)}A\right)\geq
\mu(A)^{k+1}-\varepsilon\right\}
$$
is syndetic.
\end{theorem}
We stress  that we do not assume ergodicity for this result. This
sharply contrasts the behavior of a family of linear integer
polynomials. For example when $p_i(n)=in$ for $i=1,\ldots,k$,
 it was shown  in~\cite{BHK} that the analogous result fails for
certain ergodic transformations when $k\geq 4$ and also fails for
certain nonergodic transformations when $k\geq 2$.

\section{Combinatorial Results}

Furstenberg~\cite{Fu1} established the connection between
combinatorial number theory and ergodic theory, showing that
regularity properties of subsets of integers with positive density
correspond to multiple recurrence properties of measure preserving
systems. This is reflected in what has become known as the
Correspondence Principle (first introduced in \cite{Fu1} and given
in the form below in \cite{BL}):

\begin{Correspondence1}
Let $\Lambda\subset \N$. There exist a measure preserving system
$(X,\mathcal{X},\mu, T)$ and $A\in\mathcal{X}$ such that
$\mu(A)=\bar{d}(\Lambda)$ and
$$
\bar{d}\bigl(\Lambda \cap (\Lambda+n_1)\cap\cdots\cap
(\Lambda+n_m)\bigr) \geq \mu(A\cap T^{n_1}A\cap\cdots \cap
T^{n_r}A)
$$
for all $r\in\mathbb{N}$ and all $n_1,\ldots,n_r\in\mathbb{Z}$.
\end{Correspondence1}

As an immediate corollary of Theorem~\ref{C:T-2} and Furstenberg's
Correspondence Principle, for rational independent polynomials
  we have tight lower bounds for the upper
densities in~\eqref{E:density} for every $k\in \N$. This result is
known to be false for $k\geq 4$ linear polynomials (see
~\cite{BHK}):
\begin{theorem}
\label{T:lower-bounds} Let $\Lambda\subset \N$  and  $p_1, \ldots,
p_k$ be rationally  independent integer polynomials with
$p_i(0)=0$ for $i=1,\ldots,k$. Then for every $\varepsilon>0$, the
set
\begin{equation}\label{E:bounds}
\{n\in\N\colon \quad
\bar{d}\bigl(\Lambda\cap(\Lambda+p_1(n))\cap\cdots\cap (\Lambda
+p_k(n))\bigr)\geq \bar{d}(\Lambda)^{k+1}-\varepsilon\}
\end{equation}
is syndetic.
\end{theorem}

We give an example to show that the lower bounds given in
\eqref{E:bounds} are tight.  A set $\Lambda\subset \N$ is called
normal if its indicator function ${\bf 1}_{\Lambda}$ contains
every string of zeros and ones of length $k$ with frequency
$2^{-k}$.  For any such set $\Lambda$ we have that
$$
\bar{d}\bigl(\Lambda\cap(\Lambda+n_1)\cap\cdots\cap
(\Lambda+n_k)\bigr)=\bar{d}(\Lambda)^{k+1}=1/2^{k+1}
$$
for all choices of nonzero distinct integers $n_1,\ldots,n_k$,
meaning that  \eqref{E:bounds} cannot be improved.

We remark that Furstenberg's correspondence Principle and, as a
consequence, Theorem~\ref{T:lower-bounds} hold if one replaces the
upper density $\bar{d}$ with the upper Banach density $d^*$
defined by
$d^*(\Lambda)=\lim_{N\to\infty}\sup_{M\in\mathbb{N}}|\Lambda\cap[M,M+N)|/N$
(the limit exists by subadditivity).

Szemer\'edi's Theorem has the following finite version:
given a length $k$ of a progression and density $\delta > 0$,
there
exists some $N(k,\delta)$ such that for all $N\geq N(k,\delta)$, any
subset of $\{1, \ldots, N\}$ having at least $\delta N$ elements
contains an arithmetic progression of length $k$.  In~\cite{BHK}, the
authors asked if one can strengthen this to showing that for all
$k\in\N$,
$\delta > 0$ and $\varepsilon > 0$, there exists $N(k, \varepsilon,
\delta)$ such that for all $N\geq N(k, \varepsilon, \delta)$, any
subset
of $\{1, \ldots, N\}$ with at least $\delta N$ elements contains at
least
$(1-\varepsilon)\delta^kN$ arithmetic progressions of length $k$ with
the
same common difference.  Their results show that the answer is no for $k\geq
5$
and they show that a weaker condition holds for $k=3$ and $k=4$.
Green~\cite{green} answered the (stronger) question affirmatively for
$k=3$ and $k=4$ remains open.  Given Theorem \ref{T:lower-bounds},
it is natural to ask  whether
 a similar result holds for  independent
polynomial configurations. We show that this is the case:

\begin{theorem}
Let $p_1,\ldots,p_k$ be rationally  independent integer
polynomials with $p_i(0)=0$ for $i=1,\ldots,k$. For every
$\delta>0$ and $\varepsilon>0$ there exists  $N(\varepsilon,
\delta)$, such that for all  $N>N(\varepsilon, \delta)$,  any
integer subset  $\Lambda\subset[1,N]$ with $|\Lambda|\geq\delta
N$  contains at least $(1-\varepsilon)\delta^{k+1} N$
configurations of the form $ \{x,x+p_1(n),\ldots,x+p_k(n)\} $ for
some fixed $n\in\N$.
\end{theorem}
\begin{proof}
 Suppose that
the result fails. Then there exist $\delta_0,\varepsilon_0>0$, an
integer  sequence $N_m\to\infty$, and integer subsets
$\Lambda_m\subset[1,N_m]$ such that
\begin{equation}
\label{finite1}
|\Lambda_m|\geq \delta_0N_m
\end{equation}
and
\begin{equation}
\label{finite2}|
\Lambda_m\cap(\Lambda_m+p_1(n))\cap\cdots\cap(\Lambda_m+p_k(n))|<
(1-\varepsilon_0)\delta_0^{k+1} N_m\
\end{equation}
for every $m,n\in\N$. We construct a measure preserving system
that has bad recurrence properties and then  obtain a
contradiction from Theorem~\ref{C:T-2}.

 For $m\in\N$ set $\Lambda^0_m=\Lambda^c_m$ and $\Lambda^1_m=\Lambda_m$.
  Using a
diagonal argument we can find a subsequence of $\{N_m\}_{m\in\N}$,
which for convenience we call again $\{N_m\}_{m\in\N}$, such that
the limit
$$
\lim_{m\to\infty} \frac{|(\Lambda^{i_1}_m+n_1)\cap
(\Lambda^{i_2}_{m}+n_2)\cap\cdots\cap (\Lambda^{i_r}_m+n_r)\cap
[1,N_m] |}{N_m}
$$
 exists for every $r\in\N$, $n_1,\ldots, n_r\in\bbZ$, and
$i_1,\ldots,i_r\in\{0,1\}$.

On the sequence space $(X=\{0,1\}^\bbZ,\mathcal{X})$, where
$\mathcal{X}$ is the Borel $\sigma$-algebra,  we define a measure
$\mu$ on cylinder sets as follows:
\begin{multline*}
\mu(\{x_{n_1}=i_{1},x_{n_2}=i_{2},\ldots,x_{n_r}=i_{r}\})=\\
\lim_{m\to\infty} \frac{|(\Lambda^{i_1}_m+n_1)\cap
(\Lambda^{i_2}_{m}+n_2)\cap\cdots\cap (\Lambda^{i_r}_m+n_r)\cap
[1,N_m] |}{N_m}
\end{multline*}
where $n_1, n_2, \ldots, n_r\in\bbZ$, and $i_{1}, i_2, \ldots,
i_{r}\in\{0,1\}$. The finite dimensional statistics are consistent
and so we can extend this to a probability measure on
$\mathcal{X}$ using Kolmogorov's Extension Theorem. Then the shift
transformation $T$ defined by
$$T(\{x(j)\}_{j\in\bbZ})=\{x(j+1)\}_{j\in\bbZ}
$$
preserves the measure $\mu$ and gives rise to a measure preserving
system $(X,\mathcal{X},\mu,T)$. If $ A=\{x\colon x(0)=1\}$, using
the definition of $\mu$ we see that
\begin{align*}
\mu(A\cap T^{p_1(n)}A\cap\cdots \cap T^{p_k(n)}A)&=
\mu(\{x_{0}=1,x_{p_1(n)}=1,\ldots,x_{p_k(n)}=1\})\\
\notag &=\lim_{m\to\infty} \frac{|\Lambda_m\cap
(\Lambda_{m}+p_1(n))\cap\cdots\cap (\Lambda_m+p_k(n))|}{N_m},
\end{align*}
for every $n\in\N$. Combining this with \eqref{finite1} and
\eqref{finite2} we find that
\begin{equation}\label{finite6}
\mu(A\cap T^{p_1(n)}A\cap\cdots \cap T^{p_k(n)}A))\leq
(1-\varepsilon_0)\delta_0^{k+1}\leq (1-\varepsilon_0)\mu(A)^{k+1}\
\end{equation}
for all $n\in\N$. This contradicts Theorem~\ref{C:T-2} and
completes the proof.
\end{proof}

\section{Characteristic factors and multiple recurrence result}\label{S:key}
\subsection{Preliminaries}\label{S:preliminaries}
By a measure preserving system we mean a quadruple
$(X,\mathcal{X},\mu, T)$, where $(X,\mathcal{X},\mu)$  is a
probability space and $T\colon X\to X$ is a measurable map such
that $\mu(T^{-1}A)=\mu(A)$ for all $A\in\mathcal{X}$. Without loss
of generality we can assume that the probability space is
Lebesgue.
 A \emph{factor}
of the measure preserving system $(X,\cB, \mu, T)$ can be defined
in any of the following three equivalent ways: it is a
$T$-invariant sub-$\sigma$-algebra $\mathcal{D}$ of $\cX$, it is a
$T$-invariant sub-algebra $\mathcal{F}$ of $L^\infty(X)$, or it is
a system $(Y, \cY, \nu, S)$ and a measurable map $\pi\colon X'\to
Y'$, where $X'$ is a $T$-invariant set and $Y'$ is an
$S$-invariant set of full measure, such that $\mu\circ\pi^{-1} =
\nu$ and $S\circ\pi(x) = \pi\circ T(x)$ for $x\in X'$. By setting
$\mathcal{F}=L^\infty(\mathcal{D})$, we see that the first
definition implies the second. Conversely, given $\mathcal{F}$ we
define $\cD$ to be the $\sigma$-algebra generated by
$\mathcal{F}$-measurable sets. The equivalence between the first
and third definition is seen by identifying  $\mathcal{D}$ with
$\pi^{-1}(\cY)$.  In a slight abuse of terminology, when any of
these conditions holds, we say that $Y$ (or the appropriate
$\sigma$-algebra of $\cX$) is a factor of $X$ and call $\pi\colon
X'\to Y'$ the factor map. If a factor map $\pi\colon X'\to Y'$ is
also injective, then we say that the systems $(X,\cB, \mu, T)$ and
$(Y, \cY, \nu, S)$ are \emph{isomorphic}.

If $\cY$ is a $T$-invariant sub-$\sigma$-algebra of $\cX$ and
$f\in L^2(\mu)$, we define the \emph{conditional expectation
$\mathbb{E}(f|\cY)$ of $f$ with respect to $\cY$} to be the
orthogonal projection of $f$ onto $L^2(\cY)$. We frequently
use the identities
$$
\int \mathbb{E}(f|\cY) \ d\mu= \int f\ d\mu, \quad
T\,\mathbb{E}(f|\cY)=\mathbb{E}(Tf|\cY)\ .
$$

For each $r\in\mathbb{N}$, we define $\mathcal{K}_{r}$ to be the
factor induced by the algebra
$$\{f\in L^\infty(\mu):T^rf=f\}\ .$$
We define $\mathcal{K}_{rat}$ to be the factor induced by the
algebra generated by the functions
$$
\{f\in L^\infty(\mu):T^rf=f \text{ for some } r\in\mathbb{N}\}\ .
$$
The Kronecker factor $\mathcal{K}$  is induced by the algebra spanned by the
 bounded
eigenfunctions of $T$.

The transformation $T$ is \emph{ergodic} if $\cK_1$ consists only
of constant functions and $T$ is \emph{totally ergodic} if
$\cK_{rat}$ consists only of constant functions. The von Neumann
Ergodic Theorem states that if $T$ is ergodic and $f\in L^2(\mu)$,
then
\begin{equation}\label{E:ergodic}
\lim_{N\to\infty}\frac{1}{N}\sum_{n=0}^{N-1} T^nf=\int f\ d\mu \ ,
\end{equation}
with the convergence taking place in $L^2(\mu)$.

Every measure preserving system $(X,\cX,\mu,T)$ has an
\emph{ergodic decomposition}, meaning that we can write $\mu=\int
\mu_t\ d\lambda(t)$, where  $\lambda$ is a probability measure on
$[0,1]$ and $\mu_t$ are $T$-invariant probability measures on
$(X,\cX)$ such that the systems $(X,\cX,\mu_t,T)$ are ergodic for
$t\in [0,1]$.

 If  $G$ is a $k$-step
nilpotent Lie group and $\Gamma$ is a cocompact subgroup, then
$X=G/\Gamma$ is called a $k$-\emph{step nilmanifold}. There exists
a unique probability measure $m$ on $X$ (the \emph{Haar measure})
that is invariant under left translations. If $a\in G$, then the
measure preserving system $(X,\mathcal{X},m,T_a)$ defined by the
transformation $T_a(g\Gamma)=(ag)\Gamma$ is called a
\emph{nilsystem}. Every unipotent affine transformation on a
compact abelian Lie group (with the Borel $\sigma$-algebra and
 the Haar measure) induces a system
that is isomorphic to a nilsystem, but these are not the only
examples of nilsystems.

We say that the system $(X,\cX,\mu,T)$ is an \emph{inverse limit
of a sequence of factors} $(X,\cX_j,\mu,T)$ if
$\{\cX_j\}_{i\in\mathbb{N}}$ is an increasing sequence of
$T$-invariant sub-$\sigma$-algebras such that
$\bigvee_{j\in\N}\mathcal{X}_j=\mathcal{X}$ up to sets of measure
zero. If in addition for every $j\in\N$ the factor system
$(X,\cX_j,\mu,T)$ is isomorphic to a nilsystem of order $k$, we
say that $(X,\cX,\mu,T)$ is an \emph{inverse limit of nilsystems
of order} $k$.

\subsection{Characteristic factors }
\label{sec:char}
A key ingredient in the
proof of Theorem~\ref{T:characteristic}  is the following result
of the authors:

\begin{theorem}[\cite{FrK1}]\label{T:product}
Let $(X,\cX,\mu,T)$ be an invertible  totally ergodic measure preserving
system and let $p_1,\ldots, p_k$ be rationally independent integer
polynomials. Then for $f_1, \ldots, f_k \in L^\infty(\mu)$ the
difference
\begin{equation}
\label{eq:polynomial}  \frac{1}{N-M}\sum_{ n=M}^{N-1}
T^{p_1(n)}f_1\cdot \ldots \cdot T^{p_k(n)}f_k - \prod_{i=1}^k\int
f_i\,d\mu
\end{equation}
converges to $0$ in $L^2(\mu)$ as $N-M\to\infty$.
\end{theorem}
We note that the result in~\cite{FrK1} is only stated for $M=0$,
but the same proof gives this uniform version. If $k\ge 2$ and the
polynomials $p_1,\ldots,p_k$ are not rationally independent then
there exist totally ergodic systems  and bounded functions
$f_1,\ldots,f_k$, for which the limit of the average in
\eqref{eq:polynomial} is not constant. This can be easily seen by
considering the example of an irrational rotation on the circle.

Before the proof of Theorem~\ref{T:characteristic}, we prove a
Lemma:

\begin{lemma}\label{L:spectral}
Let $(X,\cX,\mu,T)$ be a  measure preserving system with ergodic
decomposition $\mu=\int \mu_t \ d\lambda(t)$. If $f\in
L^\infty(\mu)$ satisfies $\E(f\,|\,\cK_{rat}(\mu))=0$, then
$\E(f\,|\,\cK_{rat}(\mu_t))=0$ for $\lambda$-a.e. $t$.
\end{lemma}
\begin{proof}
Let $\sigma$, $\sigma_t$  be the spectral measures of the function
$f$ with respect to the systems $(X,\cX,\mu,T)$ and
$(X,\cX,\mu_t,T)$, respectively. It is classical that
$\E(f\,|\,\cK_{rat}(\mu))=0$ if and only if $\sigma(\{r\})=0$ for
every $r\in\mathbb{Q}$. Since $\sigma=\int \sigma_t\ d\lambda(t)$,
we have that
$$
0=\sigma(\{r\})=\int \sigma_t(\{r\}) \ d\lambda(t)
$$
for every $r\in\mathbb{Q}$. Hence, for every $r\in\mathbb{Q}$ we
have $\sigma_t(\{r\})=0$ for $\lambda$-a.e. $t$. Since
$\mathbb{Q}$ is countable it follows that for $\lambda$-a.e. $t$
we have $\sigma_t(\{r\})=0$ for every $r\in\mathbb{Q}$,  and so
for $\lambda$-a.e. $t$ we have  $\E(f\,|\,\cK_{rat}(\mu_t))=0$.
\end{proof}

Every integer polynomial $p(n)$ of degree at most $d$ admits a
representation of the form $ p(n)=\sum_{i=0}^d c_i\binom{n}{i} $
for some $c_i\in\mathbb{Q}$, $i=0,\ldots,d$. Since $p(j)\in\bbZ$,
$j=0,\ldots,d$, it is immediate  that $c_i\in\mathbb{Z}$,
$i=0,\ldots,d$. A fact that we frequently use in the sequel is
that whenever $p(n)$ is an integer polynomial of degree at most
$d$, then for every $r\in\bbZ$ the polynomial $q(n)=p(d!n+r)$ has
integer coefficients. This follows easily from the aforementioned
representation.

\begin{proof}[Proof of Theorem~\ref{T:characteristic}]
 We begin with some easy reductions.
Without loss of generality we  can assume that the polynomials
$p_1,\ldots,p_k$ have integer coefficients. Indeed, suppose that
the highest degree of the polynomials $p_1,\ldots,p_k$ is $d$.
Then for every $r\in\bbZ$ the polynomial family
$\{p_i(d!n+r)\}_{i=1,\ldots, k}$ satisfies the assumptions of the
theorem and also has integer coefficients. Using the result for
$r=0,\ldots, d!-1$ and adding, we obtain the result for the family
$\{p_i(n)\}_{i=1,\ldots,k}$. Furthermore, since
$\bbE(T^jf\,|\,\cK_{rat})=T^j\bbE(f\,|\,\cK_{rat})$ for $j\in
\bbZ$, we can further assume that $p_i(0)=0$ for $i=1,\ldots,k$.

It suffices to show that if $\bbE(f_1\,|\,\cK_{rat})=0$ then the
average \eqref{E:averages} converges to zero  in $L^2(\mu)$ as
$N-M\to\infty$. If $f$ is a function with
$\E(f\,|\,\mathcal{K}_{rat})=0$ for the measure $\mu$, then by
Lemma~\ref{L:spectral}  the same property holds for almost every
measure in the ergodic decomposition of $\mu$.  Hence,  we can
assume that $T$ is ergodic.

 From \cite{Lei3} we know that a
characteristic factor for $L^2(\mu)$ convergence of the averages
\eqref{E:averages} is an inverse limit of nilsystems induced by
some $T$-invariant sub-$\sigma$-algebras  $\{\cX_j\}_{j\in\N}$.
Since $\bbE(f_1\,|\,\cK_{rat}(\cX))=0$ implies that
$\bbE(f_1\,|\,\cK_{rat}(\cX_j))=0$ for $j\in\N$,  using a standard
approximation argument we can assume that the system is a
nilsystem.

The Kronecker factor of an ergodic nilsystem is isomorphic to a
rotation on a monothetic compact abelian Lie group $G$. Every such
group has the form $\bbZ_{d_1} \times \bbT^{d_2}$ for some
positive integer $d_1$ and nonnegative integer $d_2$, where
$\bbZ_d$ denotes the cyclic group with $d$ elements. It follows
that $\cK_{rat}=\cK_{r_0}$ for some $r_0\in\N$. Hence, every
ergodic component of the transformation $T^{r_0}$ is totally
ergodic. Since $p_i(0)=0$ and $p_i$ has integer coefficients, we
have that $p_i(nr_0)=r_0q_i(n)$, where $q_i(n)$, for
$i=1,\ldots,k$,
 is again a polynomial with integer coefficients.  From
$\bbE(f_1\,|\,\cK_{rat})=0$, it follows that the function $f_1$
has integral zero on  every ergodic component of $T^{r_0}$.
Applying Theorem~\ref{T:product} on the (totally) ergodic
components of $T^{r_0}$ with the rationally independent
polynomials $q_1,\ldots, q_k$, we have that
\begin{equation}\label{E:2}
\frac{1}{N-M}\sum_{n=M}^{N-1} T^{p_1(nr_0)}f_1\cdot\ldots\cdot
T^{p_k(nr_0)}f_k
\end{equation}
converges to $0$ in $L^2(\mu)$ as $N-M\to\infty$. Moreover,
$\bbE(f_1\,|\,\cK_{rat})=0$ implies that
$\bbE(T^jf_1\,|\,\cK_{rat})=0$ for $j\in\N$ and so the limit is
zero with $p_i(nr_0+k)$ substituted for $p_i(nr_0)$ in \eqref{E:2}
for $k=0,\ldots,r_0-1$. Adding these, we have that
\eqref{E:averages} converges to $0$ in $L^2(\mu)$ as
$N-M\to\infty$.
\end{proof}

\subsection{Multiple recurrence}
\label{sec:multrec} We prove Theorem~\ref{C:T-2}.
\begin{proof}[Proof of Theorem~\ref{C:T-2}]
Suppose that the highest degree of the polynomials
$p_1,\ldots,p_k$ is $d$. Then the polynomial family
$\{p_i(d!n)\}_{i=1,\ldots, k}$ satisfies the assumptions of the
theorem and has integer coefficients. By applying the result for
this family we can assume that the polynomials $p_1,\ldots,p_k$
have integer coefficients.

Let  $\varepsilon>0$. There exists $r\in\N$ such that
\begin{equation}\label{E:approximation}
\|\bbE({\bf 1}_{A}\,|\,\cK_{r})- \bbE({\bf
1}_{A}\,|\,\cK_{rat})\|_{L^2(\mu)}\leq \frac{\varepsilon}{k+1}\ .
\end{equation}
By Theorem~\ref{T:characteristic},
\begin{align}
\label{E:lim} \lim_{N-M\to\infty}\frac{1}{N-M}\sum_{n=M}^{N-1} &
\mu(
A\cap T^{p_1(nr)}A\cap\cdots\cap T^{p_k(nr)}A)=\\
\notag \lim_{N-M\to\infty}\frac{1}{N-M}\sum_{n=M}^{N-1} & \int
\bbE({\bf 1}_{A}\,|\,\cK_{rat})\cdot T^{-p_1(nr)} \bbE({\bf
1}_{A}\,|\,\cK_{rat})\cdot \ldots \cdot T^{-p_k(nr)}\bbE({\bf
1}_{A}\,|\,\cK_{rat}) \ d\mu\ .
\end{align}
For every choice of integers $a_0,\ldots,a_k$, we have
\begin{align*}
& \Bigl| \int \prod_{i=0}^kT^{a_i}\bbE({\bf
1}_{A}\,|\,\mathcal{K}_{rat}) \ d\mu -\int
\prod_{i=0}^kT^{a_i}\bbE({\bf 1}_{A}\,|\,\mathcal{K}_{r})  \ d
\mu\Bigr|
\\
 &
\leq  \int \sum_{i=0}^k|T^{a_i}\bbE({\bf
1}_{A}\,|\,\mathcal{K}_{rat}) -
T^{a_i}\bbE({\bf 1}_{A}\,|\,\mathcal{K}_{r})|\ d \mu\\
& =  \sum_{i=0}^k \int |\bbE({\bf 1}_{A}\,|\,\mathcal{K}_{rat}) -
\bbE({\bf 1}_{A}\,|\,\mathcal{K}_{r})|\ d \mu \\
& \leq    \sum_{i=0}^k\| \bbE({\bf 1}_{A}\,|\,\mathcal{K}_{rat})-
\bbE({\bf 1}_{A}\,|\,\mathcal{K}_{r})\|_{L^2(\mu)} \leq
\varepsilon
\end{align*}
by~{(\ref{E:approximation})}.  It follows that the limit in
\eqref{E:lim} is greater than or equal to
\begin{align*}
\lim_{N-M\to\infty}\frac{1}{N-M}\sum_{n=M}^{N-1} &\int \bbE({\bf
1}_{A}\,|\,\cK_{r})\cdot T^{-p_1(nr)}\bbE({\bf
1}_{A}\,|\,\cK_{r})\cdot \ldots \cdot
T^{-p_k(nr)}\bbE({\bf 1}_{A}\,|\,\cK_{r}) \ d\mu-\varepsilon\\
=&\int \bbE({\bf 1}_{A}\,|\,\cK_{r})^{k+1} \ d\mu-\varepsilon \
\geq \mu(A)^{k+1}-\varepsilon,
\end{align*}
where the last equality holds since $r$ divides $p_i(nr)$ for
$i=1,\ldots,k$,  and every $\cK_r$-measurable function is $T^r$
invariant.
\end{proof}

\section{Correlations of independent polynomial iterates and
nilsequences}\label{SS:nil} We now prove the Structure
Theorem~\ref{T:nilsequence} for  multicorrelation sequences of
independent polynomials. We start with some definitions
from~\cite{BHK}:
\begin{definition}
Let $k\geq1$ be an integer and let $X=G/\Gamma$ be a $k$-step
nilmanifold. Suppose that  $\phi$ is a continuous complex valued
function on $X$, $a\in G$, and  $x_0\in X$. The sequence
$\{\phi(a^nx_0)\}_{n\in\mathbb{N}}$ is called a \emph{basic}
$k$-\emph{step nilsequence}. A $k$-\emph{step nilsequence} is a
uniform limit of basic $k$-step nilsequences.
\end{definition}
\begin{definition} Let $\{a_n\}_{n\in\mathbb{N}}$ be a bounded
sequence of complex numbers. We say that $a_n$ \emph{tends to zero
in uniform density}, and write $UD$-$\lim{a_n}=0$, if
$$
\lim_{N-M\to\infty}\frac{1}{N-M}\sum_{n=M}^{N-1}|a_n|=0 \ .
$$
\end{definition}

 Before the proof, we begin with a Lemma:
\begin{lemma}\label{L:UD}
Let $(X,\cX,\mu,T)$ be a  measure preserving system,
$p_1,\ldots,p_k$ be rationally independent integer polynomials,
and $f_0,f_1,\ldots,f_k\in L^{\infty}(\mu)$. Then
\begin{equation}\label{E:UD}
\text{UD-}\lim\Big(\int f_0\cdot T^{p_1(n)}f_1\cdot \ldots \cdot
T^{p_k(n)}f_k \ d\mu- \int \tilde{f}_0\cdot
T^{p_1(n)}\tilde{f}_1\cdot \ldots \cdot T^{p_k(n)}\tilde{f}_k \
d\mu\Big)=0\ ,
\end{equation}
where $\tilde{f_i}=\bbE(f_i\,|\,\cK)$, $i=0,1,\ldots,k$ and $\cK$
is the Kronecker factor of the system.
\end{lemma}
\begin{proof}
It suffices to show that if $\bbE(f_i\,|\,\cK)=0$ for some
$i\in\{1,\ldots,k\}$, then the $UD$-limit in \eqref{E:UD} is zero.
Without loss of generality, we can assume that $i=1$. We apply
Theorem~\ref{T:characteristic} to the product system induced by
$T\times T$ acting  on $X\times X$. From  \cite{Fu2} (Lemma 4.18)
we know that $f \in\cK(X\times X)$ if and only if it has the form
$$
f(x,x')=\sum_{n\in\mathbb{N}} c_n  \, g_n(x)\cdot h_n(x')
$$
where $g_n,h_n\in \cK(X)$ and $c_n\in\mathbb{C}$ for
$n\in\mathbb{N}$. Since $\bbE(f_1\,|\,\cK(X))=0$, it follows that
$\bbE(f_1\otimes
 \bar{f}_1\,|\,\cK(X\times X))=0$ which implies that $\bbE(f_1\otimes
 \bar{f}_1\,|\,\cK_{rat}(X\times X))=0$. Hence, the average
$$
\frac{1}{N-M}\sum_{n=M}^{N-1} (T\times T)^{p_1(n)}(f_1\otimes
\bar{f}_1) \cdot\ldots\cdot (T\times T)^{p_k(n)}(f_k\otimes
\bar{f}_k)
$$
converges to zero in $L^2(\mu\times\mu)$ as $N-M\to \infty$. It
follows that
$$
\lim_{N-M\to\infty} \frac{1}{N-M}\sum_{n=M}^{N-1} \Big|\int
f_0\cdot T^{p_1(n)}f_1\cdot\ldots\cdot T^{p_k(n)}f_k\
d\mu\Big|^2=0\ .
$$
and this completes the proof.
\end{proof}
\begin{proof}[Proof of Theorem~\ref{T:nilsequence}]
 By Lemma~\ref{L:UD}, we can assume that $\cX=\cK$. Since the system is
 ergodic and coincides with its Kronecker factor
 we can assume that $T$ is a rotation on a compact abelian
group $G$. Every compact abelian group is an inverse limit of
compact abelian Lie groups and so using an easy approximation
argument, such as the one used in \cite{BHK} (see page 296), we can
further assume that $G$ is Lie.

Suppose now that $G$ is a compact abelian Lie group with Haar
measure $m$ and that $T\colon G\to G$ is given by $T(g)=g+a$ for
some $a\in G$.  For $j=0,\ldots,d$ we have that
$p_j(n)=\sum_{i=0}^dc_{i,j}\ \binom{n}{i}$ for some
$c_{i,j}\in\bbZ$. We construct the advertised transformation
$S\colon G^d\to G^d$ and the continuous function $\phi\colon
G^d\to \mathbb{C}$ as follows: $S$ is defined by
$$
S\big(g_1,g_2,\ldots,g_d\big)=\big(g_1+a,g_2+g_1,\ldots,
g_d+g_{d-1}\big)\ ,
$$
and the continuous function $\phi$ is defined by
 $$
  \phi(g_1,\ldots,g_d)=\int f_0(g)\cdot
\prod_{i=1}^k f_i(g+c_{i,0}a+\sum_{j=1}^d c_{i,j}g_i) \ dm(g)\ .
$$
Note that $S$ is unipotent since all its eigenvalues are $1$.
 It is easy to check that
$$
S^n(0,\ldots,0)=\Big(\binom{n}{1}a,\ldots,\binom{n}{d}a\Big)\ ,
$$
and so
\begin{align*}
\phi\bigl(S^n(0,\ldots,0)\bigr)=&\phi\Big(\binom{n}{1}a,\ldots,\binom{n}{d}a\Big)\\
=&\int f_0(g)\cdot f_1(g+p_1(n)a)\cdot \ldots\cdot
f_k(g+p_k(n)a) \ dm(g)\\
=&\int f_0\cdot T^{p_1(n)}f_1\cdot \ldots \cdot T^{p_k(n)}f_k \
dm= a_n\ .
\end{align*}
The system $(G^d,S)$  is topologically conjugate to a $d$-step
nilsystem, meaning that there exist a $d$-step nilmanifold
$H/\Gamma$, an $a\in H$, and an invertible continuous map
$\pi\colon G^d\to H/\Gamma$ such that $S=\pi^{-1}\circ T_a\circ
\pi$, where $T_a$ is defined by $T_a(g\Gamma)=(ag)\Gamma$. It
follows that
$$
\phi\big(S^n(0,\ldots,0)\big)=\phi'(T_a^nx_0)
$$
where $\phi'=\phi\circ \pi^{-1}$ is a continuous function on
$H/\Gamma$ and $x_0=\pi(0,\ldots,0)\in H/\Gamma$. This completes
the proof.
\end{proof}
We illustrate the construction of this  proof with an example:
\begin{example}
Suppose that $k=2$ and $p_1(n)=2n+1$, $p_2(n)=n^2/2-n/2$,
$G=\bbT$, and $T\colon \bbT\to \bbT$ is given by $T(t)=t+\alpha
\pmod{1}$ for some irrational $\alpha\in\bbT$. Then
$$
a_n=\int f_0\big(t\big) \cdot f_1\big(t+(2n+1)\alpha\big) \cdot
f_2\Big(t+\big(\binom{n}{2}-\binom{n}{1}\big)\alpha\Big)\ dt\ ,
$$
$S\colon \bbT^2\to\bbT^2$ is defined by
$$
S(t_1,t_2)=(t_1+\alpha,t_2+t_1)\ ,
$$
and $\phi\colon \bbT^2\to \mathbb{C}$ is defined by
$$
\phi(t_1,t_2)=\int f_0(t) \cdot f_1(t+\alpha+2t_1) \cdot
f_2(t-t_1+t_2)\ dt\ .
$$

\end{example}


\end{document}